\documentclass[12pt]{article}
\usepackage{graphicx}
\usepackage{amsfonts}

\font\tendb=msbm10 at 12pt \font\sevendb=msbm10 at 9pt
\font\fivedb=msbm10 at 7pt
\newfam\dbfam
\textfont\dbfam=\tendb \scriptfont\dbfam=\sevendb
\scriptscriptfont\dbfam=\fivedb
\def\db{\fam\dbfam\tendb}
\font\eufm=eufm10\font\eufms=eufm10\font\eufmss=eufm10\newfam\eufam
\textfont\eufam=\eufm\scriptfont\eufam=\eufms\scriptscriptfont\eufam=\eufmss

\font\tendbb=msbm10 at 12pt \font\sevendbb=msbm7 at 9pt
\font\fivedbb=msbm5 at 6pt
\newfam\dbbfam
\textfont\dbbfam=\tendbb \scriptfont\dbbfam=\sevendbb
\scriptscriptfont\dbbfam=\fivedbb

 \def \Z {{\db Z}}
\def \F {{\db F}}

 \newcommand{\func}[5]{
\begin{array}{cccl}
#1:& #2 & \longrightarrow & #3 \\
   & #4& \longmapsto & #5
\end{array}}
\newcommand{\doublefunc}[7]{
\begin{array}{cccl}
#1:& #2 & \longrightarrow & #3 \\
   & #4& \longmapsto & #5\\
   &#6& \longmapsto & #7
\end{array}}

\def \fin {\hfill \framebox(7,7) }
 
\font\tenMmm=eusm10 at 12pt
\newfam\Mmmfam
\textfont\Mmmfam=\tenMmm

\def\illu #1 by #2 (#3){
  \vbox to #2{
    \hrule width #1 height 0pt depth 0pt
    \vfill
    \special{illustration #3} 
    }
  }

\begin{document}
\null \vspace{2cm}
\begin{center}
{\large {\bf Equivariant Khovanov homology associated with symmetric links}}\\
 Nafaa Chbili
 \begin{footnotesize}\footnote{This work was partially done while I was visiting
 KAIST supported by a fellowship from the project  BK21.
  I would like to thank Professor K. Hyoung Ko for his
  kind hospitality.}\\
 Department of Mathematics\\
 Korea Advanced Institute of Science and  Technology\\
 Daejeon, 305-701, Korea\\
 E-mail: chbili@knot.kaist.ac.kr
 \end{footnotesize}
 \end{center}
\begin{abstract}
Let $\Delta$ be a trivial knot in the three-sphere. For every finite
cyclic group $G$ of odd order, we construct a $G$-equivariant
Khovanov homology  with coefficients in the filed $\F_{2}$. This
homology is an invariant of links up to isotopy in $(S^{3},\Delta)$.
Another interpretation is given using the categorification of the
Kauffman bracket skein module of the solid torus. Our techniques
apply  in the case of  graphs as well to define an equivariant
version of the
graph homology which categorifies  the chromatic polynomial.\\
\emph{Key words.} Khovanov homology,  group action,
equivariant Jones polynomial, skein modules.\\
\emph{MSC.} 57M25.\\
\end{abstract}
\subsection*{1- Introduction}
In the late nineties, M. Khovanov \cite{Kh} introduced an invariant
of isotopy classes of oriented  links in the three-sphere, now
widely known as the Khovanov homology. This invariant takes the form
of bigraded homology groups whose polynomial Euler  characteristic
is the Jones polynomial. Namely, if $L$ is an oriented link and
$H^{*,*}(L)$ is its Khovanov homology  with integral coefficients,
then the Jones polynomial of $L$ is given by the following formula
 $$V(L)(q)=\displaystyle\sum_{i,j}(-1)^iq^{j}\mbox{rank}H^{i,j}(L),$$
 where  $V(L)(q)$ is  the augmented version of the Jones polynomial,
 equal to $(q+q^{-1})$ times the
  original Jones invariant
\cite{Jo}. The original definition of Khovanov homology is
complicated and overloaded with algebraic details.  Viro \cite{Vi}
suggests  an elementary combinatorial approach to define the
Khovanov homology. This approach has proved to be useful  in several
works. For instance, it was used in \cite{APS} to construct an
homology theory for framed links in $I-$bundles over surfaces.
This theory  categorifies the Kauffamn bracket skein module \cite{Pr1}.\\
Quantum invariants of links have proved to be a powerful tool in the
study of the symmetry of links. For instance, the Jones and the
HOMFLY  polynomials satisfy certain necessary conditions which
helped determine the  symmetries of some links \cite{Mu,Tr,Pr2}. Our
main goal in
this paper is to investigate the behavior of the   Khovanov homology of links
with $\Z/p\Z-$symmetry.\\
Let $\Delta$ be a trivial knot in $S^3$ and let $L$ be a link in
$S^3$ such that $L$ does not intersect $\Delta$. Let $\tilde L$ be
the covering link of $L$ in the  $p-$fold cyclic  cover branched
over $\Delta$. Obviously, the group $G=\Z/p\Z$ acts on $(S^3,\tilde
L)$. Let   $D$ be a diagram of the link $L$ and let $\tilde D$ be a
symmetric diagram of $\tilde L$. We prove that the action of $G$
extends naturally to the Khovanov chain complex of $\tilde D$ with
coefficients in $\F_2$. The homology of the quotient chain complex
is called here the \emph{$G-$equivariant Khovanov homology} of $D$,
we shall denote it by $H_G^{*,*}(D)$. Throughout this paper, two
links in $(S^3, \Delta)$ are isotopic if they are related by an
isotopy of $S^3$ keeping the knot $\Delta$ fixed. If  $E$ is a
vector space on which the group $G$ acts, then we set  $E^G$ to be
the subset
of fixed points under this action.\\

\textbf{Theorem 1.} \emph{If the order of $G$ is odd, then the
$G-$equivariant Khovanov homology $H_G^{*,*}$ is an invariant of
ambiant isotopy of oriented links in $(S^3,\Delta)$. In addition, $H_G^{*,*}(L)$ is isomorphic to
the subspace of fixed points   $H^{*,*}(\tilde L)^G$. } \\

The polynomial Euler characteristic of $H_G^{*,*}$ is an
 invariant  of ambient isotopy of   links in $S^3$ which we  call here  the \emph{$G-$equivariant Jones
 polynomial:}
 $$V_{G}(L)(q)=\displaystyle\sum_{i,j}(-1)^iq^{j}\mbox{dim}H_G^{i,j}(L).$$
 \textbf{Corollary 1.} \emph{If $V_{G}(L)\neq V(\tilde L)$, then the action
 of  $G$ in homology is not trivial}.\\
 In \cite{APS}, Asaeda, Przytycki and Sikora constructed an homology theory
 which categorifies the Kauffman bracket skein modules of $I-$bundles
 over surfaces. In the case of the solid torus, this homology
  is an invariant of framed links which associates to each framed
  link $L$ homology groups $H^{*,*,*}(L)$, where scripts are
  integers. Let $L$ be a framed link in the solid torus $S^1 \times I \times I$ and let $D$ be a diagram of $L$ in the annulus.
  Let $\tilde L$ be the pre-image of $L$  in the $p-$fold cyclic cover of the solid
  torus. Let $\tilde D$ be a symmetric diagram of $\tilde L $ in the annulus.
   We prove that the finite cyclic group  $G=\Z/p\Z$ acts on the chain  complex $(C^{*,*,*}(\tilde D),
   d)$, where coefficients are taken in $\F_2$. Thus we construct a $G-$equivariant Khovanov homology $H_G^{*,*,*}(D)$ and
we prove that this homology defines an invariant of framed links.\\
\textbf{Theorem 2.} \emph{If the order of $G$ is odd, then the
$G-$equivariant Khovanov homology $H_G^{*,*,*}$ is an invariant of
framed  links in the solid torus.} \\
\textbf{Examples.} Computing  the Khovanov homology of a link  is
not an easy task in general. The computation of the equivariant
Khovanov
homology is even more difficult. We give here some easy examples with $G=\Z/3\Z$.\\
If $L$ is a trivial knot such $\Delta\cup L$ is a trivial link, then
the only non trivial homology spaces are
$H_G^{0,3}(L)=H_G^{0,-3}(L)=\F_2$ and
$H_G^{0,1}(L)=H_G^{0,-1}(L)=\F_2$. Since  $H^{0,3}(\tilde
L)=H^{0,-3}(\tilde L)=\F_2$, $H^{0,1}(\tilde L)=H^{0,-1}(\tilde
L)=(\F_2)^3$, then  the equivariant homology of $L$ is different
from the Khovanov homology of $\tilde L$. The equivariant Jones
polynomial of $L$ is different from the Jones polynomial of $\tilde
L$, as we have $V_G(L)=q^3+q+q^{-1}+q^{-3}\neq V(\tilde
L)=q^3+3q+3q^{-1}+q^{-3}$. In conclusion, the Khovanov homology of
$L$, the Khovanov homology of $\tilde L$  and the $G-$equivariant
Khovanov homology of $L$ are
different.\\
Now, we consider the knot $L$ depicted by the picture below, where
the linking number of  $\Delta$ and  $L$ is equal to 2. The covering
link $\tilde L$ is the trefoil knot. Computations show that the
$G-$equivariant Khovanov homology of $L$ is equal to  the Khovanov
homology of $\tilde L$, the non-trivial homology spaces are listed
below:
$H_G^{0,1}(L)=H_G^{0,3}(L)=H_G^{2,5}(L)=H_G^{2,7}(L)=H_G^{3,7}(L)=H_G^{3,9}(L)=\F_2$.
The equivariant Jones polynomial of $L$ is  $V_G(L)=-q^9+q^5+q^3+q$
which is equal to the Jones polynomial of $\tilde L$.\\
\begin{center}
\includegraphics[width=3cm,height=1.5cm]{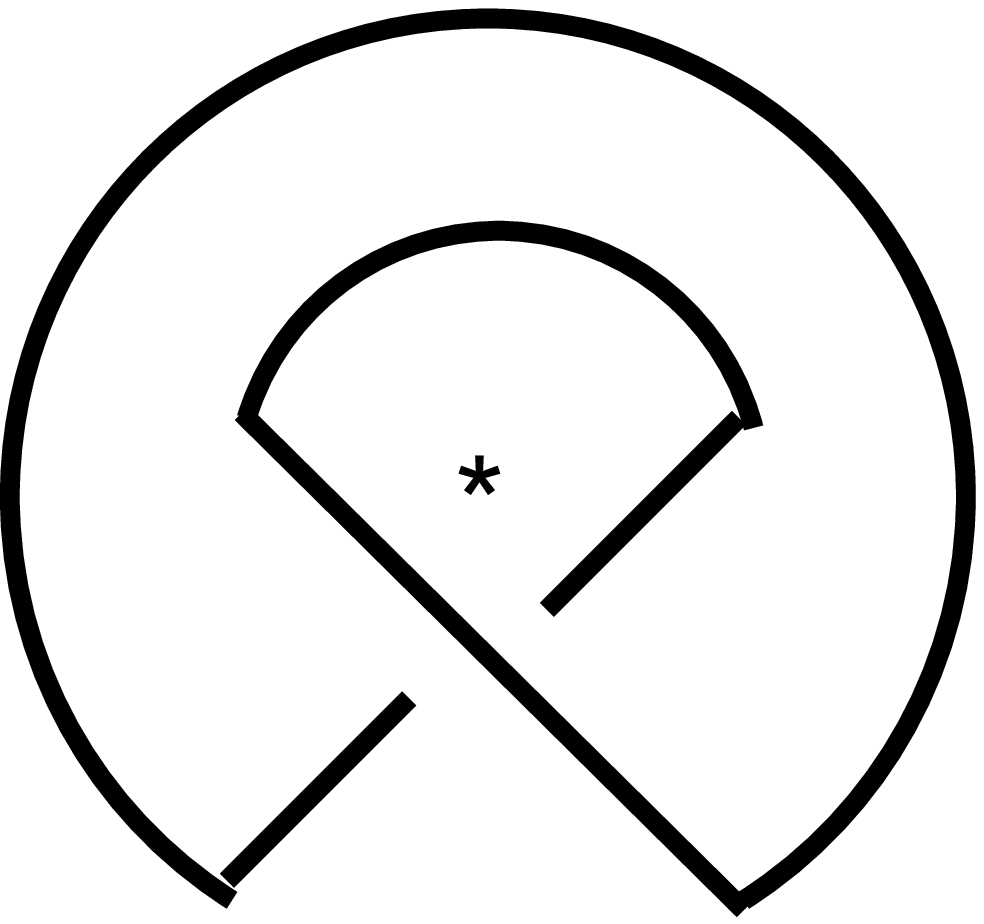}
\end{center}
Here is an outline of our paper. In Section 2, we review the
construction of Khovanov homology following \cite{Vi}. Section 3
discusses the Khavanov homology of symmetric links. In Section 4, we
review some basic properties of the transfer map in homology needed
in the sequel. The proof of  Theorem 1 is given in Section 5.
Section 6 and 7 discusses extension of our construction to framed
links in the solid torus and to graph homology.
\subsection*{2- Khovanov homology} This section is to review the
definition  of the Khovanov homology of links following the
elementary combinatorial construction introduced by Viro \cite{Vi}.
Note that coefficients will be always  in $\F_2$ and will usually be
dropped from the notation except when desired for stress.\\
Let $D$ be a link diagram with $n$ crossings. A \emph{Kauffman
state} of $D$ is an assignment of $+1$ marker or $-1$ marker to each
crossing of $D$. In a Kauffman state  the crossings of $D$ are
smoothed according to the following convention
\begin{center}
\includegraphics[width=10cm]{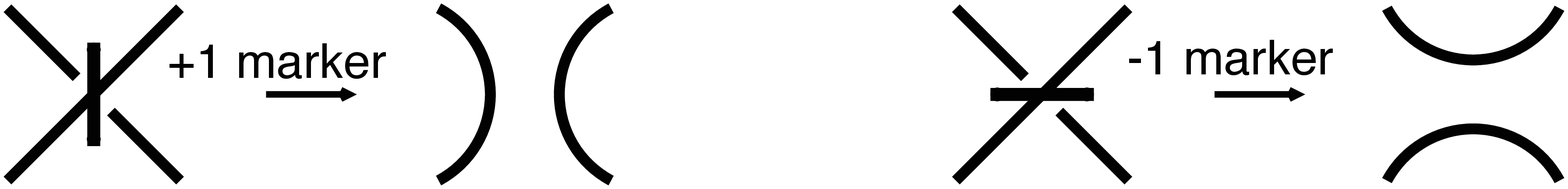}
\end{center}
\begin{center} Figure 1 \end{center}
to obtain a collection of  circles $D_s$. Let $|s|$  be the number
of circles in $D_s$ and let
 \begin{center}$\sigma(s)= \sharp \{\mbox{$+$1 markers}\}$ $-$ $\sharp
\{\mbox{$-$1 markers}\}.$ \end{center} The augmented version of
Kauffman bracket  of $D$ is the Laurent polynomial in the
indeterminate $A$  given by the following formula:
$$\prec D\succ (A)= \displaystyle\sum_{\rm{ states } \;s \;\rm{ of }\; D}
(-A)^{\sigma(s)}(-A^2-A^{-2})^{|s|}. $$
 An \emph{enhanced Kauffman state}
$S$ of $D$ is a Kauffman state $s$ with an assignment of a $+$ or
$-$ sign to each circle in $D_s$. We set $\tau(S)$ to be the
algebraic sum of signs associated to the circles
of  $D_s$.\\
If $D$ is given an orientation, then   let $w(D)$ stands for the
writhe of $D$. Now, we define:
$$\begin{array}{lll}
i(S)=&\displaystyle\frac{w(D)-\sigma(s)}{2}& \mbox{ and}\\
&&\\
 j(S)=&\displaystyle\frac{3w(D)-\sigma(s)+2\tau(S)}{2}.&
\end{array}
$$
One may check easily that both $i(S)$  and $j(S)$ are integers. Let
$i$ and $j$ be two integers, we define ${\cal S}_D^{i,j}$ to be the
set of  states of $D$ with $i(S)=i$ and $j(S)=j$. The Khovanov chain
space $C^{i,j}(D)$ is defined to be the vector space over $\F_{2}$
having
${\cal S}_D^{i,j}$ as a basis.\\
It remains to define the differential. Assume that $v$ is a crossing
of $D$, we define the partial differential $d_v$ as follows
\begin{center}
$\func{d_v^{i,j}}{C^{i,j}(D)}{C^{i+1,j}(D)}{S}{\displaystyle\sum_{\mbox{\tiny
{All states  S'}}}(S:S')_vS'}$
\end{center}
where $(S:S')_v$ is
\begin{itemize}
\item 1 if $S$ and $S'$ differ only at the crossing $v$,
where $S$ has a  $+1$ marker, $S'$ has a $-1$ marker, all the common
circles in $D_S$ and $D_{S'}$ have the same signs and in a
neighborhood of $v$, $S$ and $S'$ are as in figure 2,
\item $(S:S')$ is zero otherwise.
\end{itemize}
\begin{center}
\includegraphics[width=10cm]{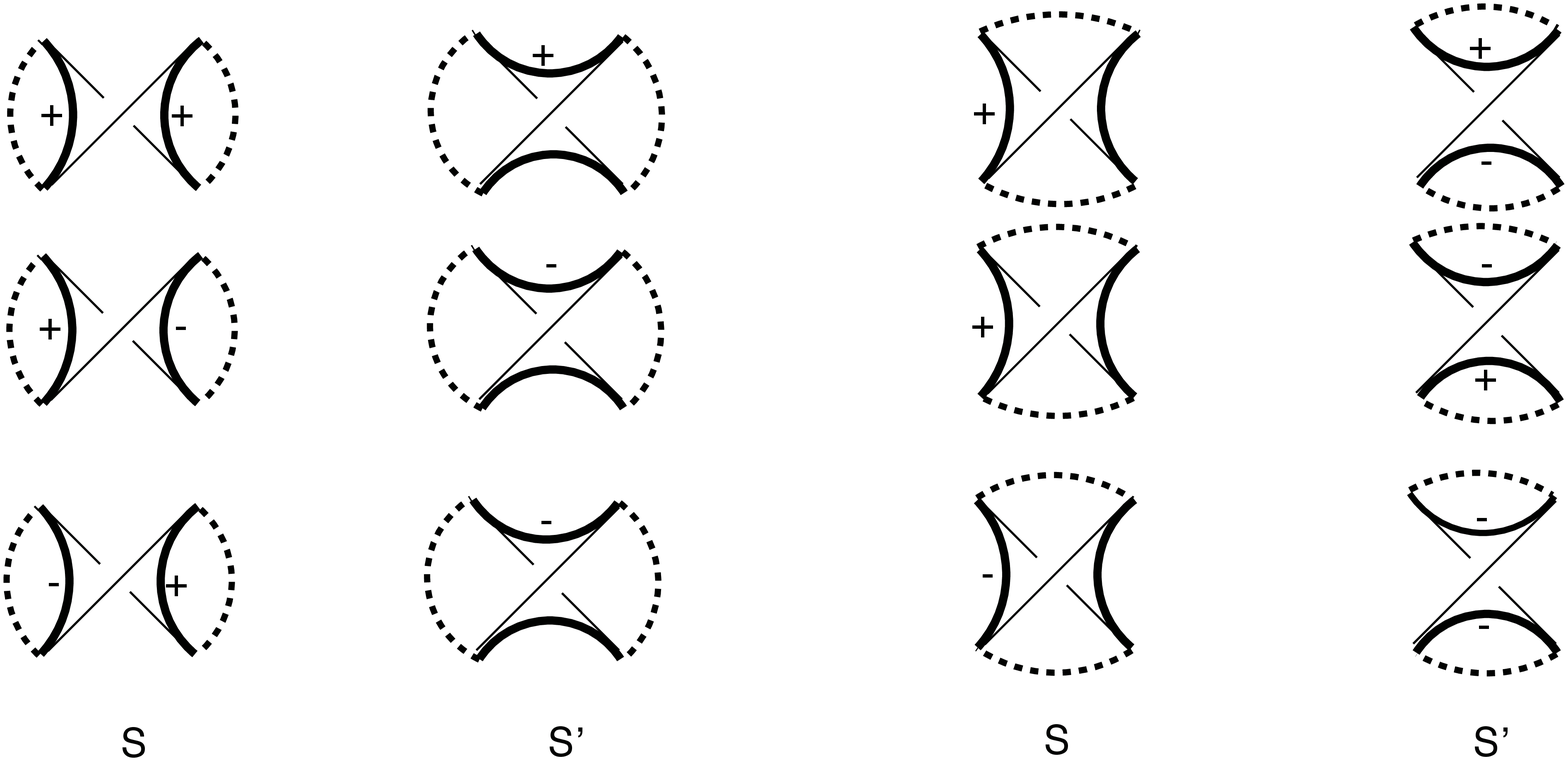}
\end{center}
\begin{center} Figure 2. \end{center}
The differential $d$ is defined by:
\begin{center}
$\func{d^{i,j}}{C^{i,j}(D)}{C^{i+1,j}(D)}{S}{\displaystyle\sum_{\mbox{\tiny
{$v$}}}d_v^{i,j}(S)}$
\end{center}
The homology $H^{*,*}(D)$ of the chain complex
$(C^{*,*}(D),d^{*,*})$ is called the Khovanov homology of $D$. This
homology is conserved under Reidemeister moves. Hence, it is  an
invariant of ambiant isotopy of links. If $L$ is an oriented link in
$S^3$, then we denote its Khovanov homology by $H^{*,*}(L)$. As we
have mentioned in the introduction, the Jones polynomial of $L$ is
obtained
as the polynomial Euler characteristic of $H^{*,*}(L)$.\\
\textbf{Framed Khovanov homology.} As it is the case for the Jones
polynomial. It is sometimes more convenient to work with framed
links when studying Khovanov homology. Viro \cite{Vi} showed that
one may define a Khovanov homology which categorifies the Kauffamn
bracket polynomial. Let $D$ be a nonoriented link diagram. With
respect to the notations of section 2, we set $p(S)=\tau(S)$ and
$q(S)=\sigma(S)-2\tau(S)$. Let $C_{p,q}(D)$ be the vector space
generated by all enhanced states with $p(S)=p$ and $q(S)=q$. We have
a chain complex $(C_{*,*}(D),d)$, where $d: C_{p,q}(D)\longmapsto
C_{p-1,q}(D)$ is defined as in the previous paragraph. If $D$ is
oriented, then we get the Khovanov homology of $C^{*,*}(D)$ by
shifting the degrees in the homology of $C_{*,*}(D)$. The advantage
of this framed version of Khovanov homology is that there is a short
exact sequence which categorifies the Kauffman bracket skein
relation \cite{Vi}. Let $D$, $D_{0}$ and $ D_{\infty}$ be three link
diagrams which are identical except in a small disk
where they are like in the following picture\\
\begin{center}
\includegraphics[width=2cm,height=1.5cm]{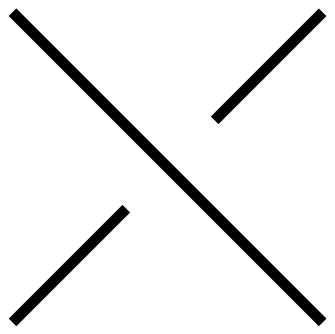} \hspace{1cm}
\includegraphics[width=2cm,height=1.5cm]{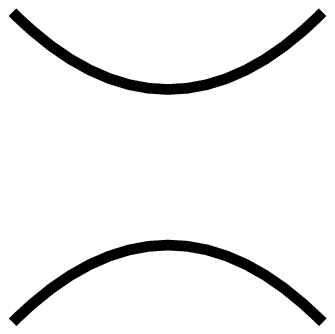} \hspace{1cm}
\includegraphics[width=2cm,height=1.5cm]{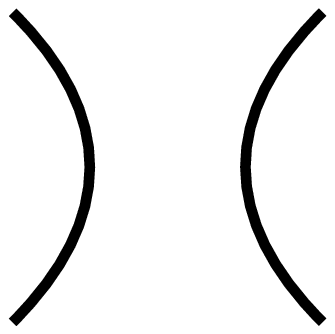}
\end{center}
\begin{center} {\sc Figure 3.} \end{center}
The following short  sequence is exact:
$$0 \longmapsto
C_{p,q}(\includegraphics[width=0.5cm,height=0.5cm]{kinfini})
\stackrel{\alpha}{\longmapsto}
C_{p,q-1}(\includegraphics[width=0.5cm,height=0.5cm]{kplus})\stackrel{\beta}{\longmapsto}
 C_{p,q-2}(\includegraphics[width=0.5cm,height=0.5cm]{kzero}) \longmapsto 0
 $$
where $\alpha$ is the chain map defined by: \includegraphics[width=0.5cm,height=0.5cm]{kinfini}$\;\longmapsto\;$\includegraphics[width=0.5cm,height=0.5cm]{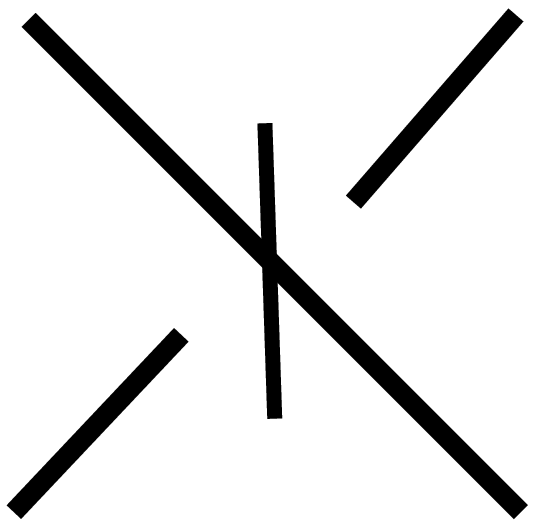}\\
and $\beta$ is defined by the following correspondence
\begin{center}
\includegraphics[width=0.5cm,height=0.5cm]{negativemarker}$\;\longmapsto\;$
0 \end{center}
\begin{center}
\includegraphics[width=0.5cm,height=0.5cm]{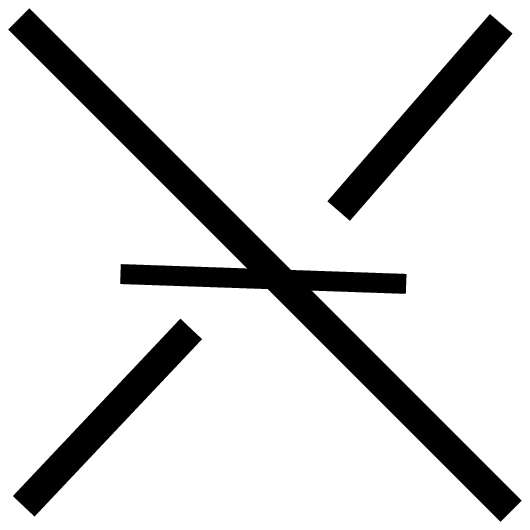}$\;\longmapsto\;$
\includegraphics[width=0.5cm,height=0.5cm]{kzero}.
\end{center}
\subsection*{3- Symmetric links and equivariant Khovanov homology} This
section is concerned with the natural question of whether the
Khovanov homology reflects the symmetry of links. In other words,
does the invariance of a link by some finite group action on the
three-sphere induce some group action on the Khovanov homology of
the link?\\
 A link $L$ in $S^3$ is said to be \emph{$p-$periodic} if
and only if there exists an orientation preserving diffeomorphism
$\varphi$ of $S^3$ such that $\varphi$ is of order $p$, the set of
fixed points of $\varphi$ is a knot disjoint from $L$ and
$\varphi(L)=L$. By the positive solution of the smith conjecture, we
may assume without loss of generality  that $\varphi$ is a rotation
of $2\pi/p-$angle around a trivial knot. Consequently, a
$p-$periodic knot admits a planar diagram which is invariant by a
planar rotation
of the same angle.\\
Let $\Delta$ be a trivial knot and let $L$ be a link in the
three-sphere such that $L\bigcap\Delta=\emptyset$. Let $\tilde L$ be
the covering link of $L$ in the $p-$fold cyclic  cover branched over
$\Delta$. Let  $\tilde D$ be a diagram of $\tilde L$ which is
invariant by a planar rotation. Such a diagram exists since the link
$\tilde L$ is $p-$periodic. Let $D$ be the quotient diagram of
$\tilde D$ under the action of the group $G=\Z/p\Z=<\varphi>$. The
rest of this section is devoted to study the Khovanov chain complex
$(C^{*,*}(\tilde D),d^{*,*})$.\\
One can see easily that the action of the rotation on the diagram
$\tilde D$ extends naturally to an action of the cyclic group $G$ on
the set of enhanced Kauffman   states of $\tilde D$. In addition,
for all enhanced state $S$ we have:
$$
 i(\varphi^{k}(S))=i(S) \; \mbox{and}\;  j(\varphi^{k}(S))=j(S)
\;\; \mbox{ for all $1\leq k \leq p$}.
$$
Consequently, the group $G$ acts on the set ${\cal S}_{\tilde
D}^{i,j}$. Since ${\cal S}_{\tilde D}^{i,j}$ is a basis for
$C^{*,*}(\tilde D)$, then this action extends naturally to an action
of  $G$ on $C^{*,*}(\tilde D)$. It remains now to check if the
action of $G$ commutes with the
differential.\\
\textbf{Lemma 3.1.} \emph{We have: $\varphi\circ d=d\circ \varphi$.}\\
\emph{Proof:} Let $S$ be an enhanced state, one can easily see that
for every crossing in $\tilde D$ we have
$(S:S')_v=(\varphi(S),\varphi(S'))_{\varphi(v)}$. Thus:
$$\begin{array}{ll}
 \varphi(d_v
(S))&=\varphi(\displaystyle\sum_{\mbox{\tiny {All states
S'}}}(S:S')_vS') \\
&=\displaystyle\sum_{\mbox{\tiny {All states
S'}}}(S:S')_v\varphi(S')\\
&= \displaystyle\sum_{\mbox{\tiny {All states
S'}}}(\varphi(S):\varphi(S'))_{\varphi(v)}\varphi(S')\\
&= \displaystyle\sum_{\mbox{\tiny {All states
T}}}(\varphi(S):T)_{\varphi(v)}T\\
&= d_{\varphi(v)} \circ \varphi(S).
\end{array}
$$
Finally, we get $\displaystyle\sum_{\mbox{\tiny {crossings
v}}}\varphi(d_v^{i,j}(S))=\displaystyle\sum_{\mbox{\tiny {crossings
v}}}d_v^{i,j}(\varphi(S))$ which means that $\varphi\circ d=d\circ
\varphi$. This ends the proof of the lemma. \fin \\

 Let
$(\overline{C^{*,*}(\tilde D)},\overline d)$ be the quotient chain
complex of $(C^{*,*}(\tilde D), d)$ by the action of $G$. The
homology of the quotient chain complex $(\overline{C^{*,*}(\tilde
D)},\overline d)$ is called the \emph{$G-$equivariant homology} of
$D$.
We denote this homology  by $H_{G}^{*,*}(D)$.\\
\textbf{Remark 3.1.} Similarly, a framed version of equivariant
Khovanov homology can be defined for non oriented diagrams. We shall
denote
it by here $H_{*,*}^G $.\\
\textbf{Remark 3.2.} If we consider Khovanov homology with
coefficients in $\Z$ as in the original definition \cite{Kh}. We
still have an action of the group $G$ on the Khovanov chain groups
but that action does not commute with the differential. This is due
to the signs that appear in the definition of the differential.
Actually, this is the raison for which we choose to work  with
coefficients in $\F_2$.
\subsection*{4- The transfer in homology} In this section, we review some
properties of the transfer map. Let $G=<\varphi>$ be the finite
cyclic  group of order $p$ and let $(C^{*},d)$ be a chain complex
with coefficients in some field $F$. Assume that $G$ acts on the
chain complex  $(C^{*},d)$ and set $(\overline {C^{*}},\overline d)$
to be the quotient chain complex. We denote by $\pi$ the canonical
surjection with respect to the action of $G$. Let $t$ be the map
from  $(\overline {C^{*}},\overline d)$ to $(C^{*},d)$ defined by
$t(\bar S)= S+\varphi(S)+ ...+\varphi^{p-1}(S)$. The map $t$ induces
a map $t_*$ from the homology of $(\overline {C^{*}},\overline d)$
to the homology of $(C^{*},d)$. This map  called the transfer has
been useful in the study of homological properties of topological
transformation groups. The following properties are extracted from
\cite{Br}.\\
\textbf{ Theorem 4.1.} \emph{The composition $\pi_*t_*$ is the
multiplication by $p$. It is an isomorphism if the field  $F$ is of
characteristic zero or prime to $p$.}\\
Obviously, the action of $G$ on the chain complex $(C^{*},d)$
induces an action of $G$ on the homology. We have:\\
 \textbf{
Theorem 4.2.} If the field $F$ is of characteristic zero or
prime to $p$, then:\\
 $\begin{array}{ccccl} &\pi_*:& H((C^{*},d))^G &
\longrightarrow & H((\overline {C^{*}},\overline d))
\end{array}$
is an isomorphism, as is\\
 $\begin{array}{ccccl}
  &t_*:&
H((\overline {C^{*}},\overline d)) \longrightarrow & H((C^{*},d))^G.
\end{array}$

\subsection*{5- Proof of Theorem 1} In this section, we shall prove that
if the order of $G$ is odd, then the $G-$equivariant Khovanov
homology does not change under Reidemeister moves $R1, R2$ and $R3$.
Note that as we consider isotopy in $(S^3,\Delta)$, then we should
consider only Reidemeister moves which are performed in a a three
ball  which does not intersect $\Delta$.
\subsubsection*{5.1- Invariance under first
Reidemeister move} Let $D$ and $D'$ be two link diagrams which are
related by a
Reidemeister move $R1$. Assume that $D$ is the diagram in the middle of figure 4 and that $D'$ is the right twisted diagram.\\
\begin{center}
\includegraphics[width=3cm,height=1cm]{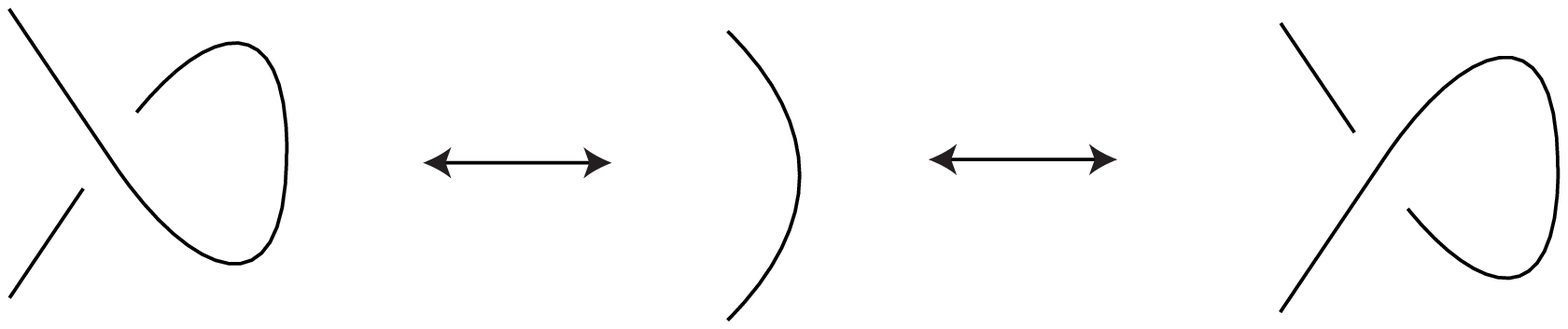}
\end{center}
\begin{center} Figure 4. \end{center}
Let $\tilde D$ and $\tilde {D'}$ be the two covering diagrams.
Obviously these two diagrams  differ by $p$ Reidemeister moves of
type 1
performed along an orbit of the action of $G$.\\
Following Viro \cite{Vi}, if two diagrams differ by a Reidemeister
move  $R1$, there is a chain map $h_v$ ($v$ is the crossing which
appears in $D'$ but not in $D$) between the two complexes which
induces an isomorphism in homology. This map is defined by
\begin{center}
\includegraphics[width=4cm,height=1cm]{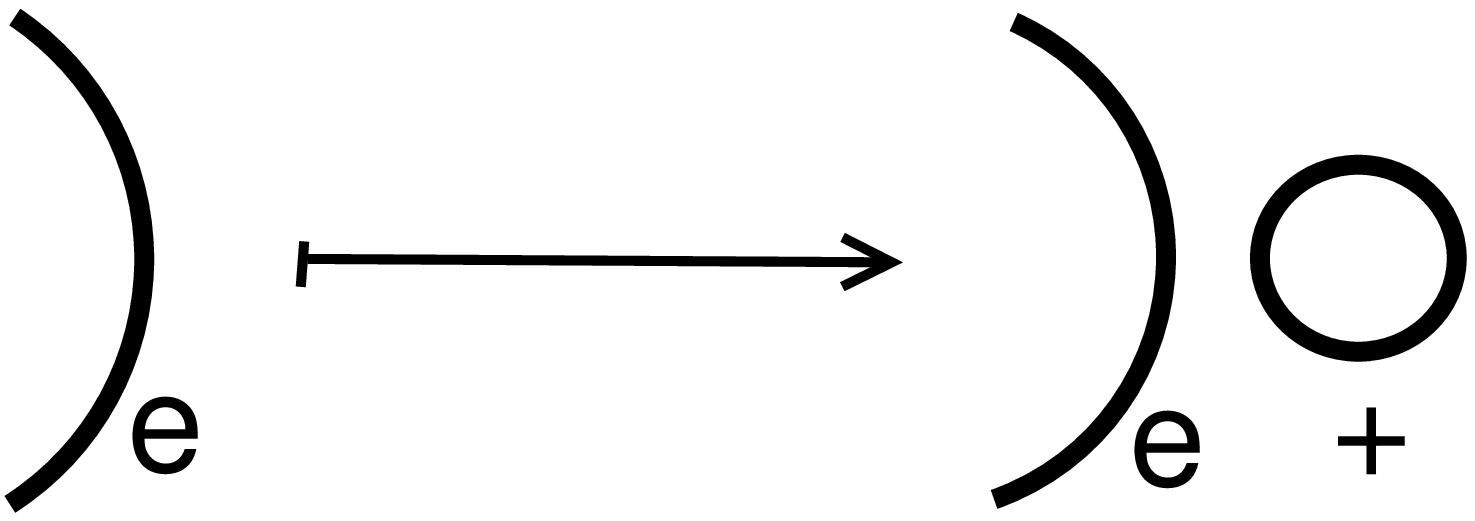}
\end{center}
The two diagrams $\tilde D$ and $\tilde {D'}$ differ by $p$
Reidemeister moves of type  R1. Let us label the crossings which
appear in $\tilde {D'}$ and not in $D$ by $v,\varphi(v),\dots
\varphi^{p-1}(v)$. Now, let $h$ be the chain map
$h_v\circ h_{\varphi(v)}\circ\dots\circ h_{\varphi^{p-1}(v)}$.\\
\textbf{Lemma 5.1.} \emph{The linear map  $h: C^{*,*}(\tilde D)
\longmapsto C^{*,*}(\tilde D') $ induces an isomorphism
in homology. In addition $h$  is $G-$equivariant.}\\
\emph{Proof.} The induced map $h_*$ is an isomorphism in homology
because it is the composition of isomorphisms. It is $G-$equivariant
due to the two elementary facts: $h_w\circ h_{w'} =h_{w'}\circ h_w$
and $\varphi \circ h_w=h_{\varphi(w)}\circ \varphi$. \fin \\

According to Lemma 5.1, the map $h$ induces a map $\overline h$ from
$(\overline{C^{*,*}(\tilde D)},\overline d)$ to
$(\overline{C^{*,*}(\tilde D')},\overline d)$. We are going to prove
that this map induces an isomorphism in homology. Note that we have
a commutative diagram \\
\begin{center}
$ {\begin{array}{ccc}
C^{*,*}(\tilde D)&\stackrel{h}{\longrightarrow}& C^{*,*}(\tilde {D'})\\
\Big\downarrow\vcenter{%
\rlap {$\pi$}}&& \Big\downarrow\vcenter{%
\rlap {$\pi'$}}\\
\overline{C^{*,*}(\tilde
D)}&\stackrel{\overline{h}}{\longrightarrow}&\overline{
C^{*,*}(\tilde {D'})}
\end{array}}$
\end{center}
which induces a commutative diagram in homology\\
\begin{center}
${\begin{array}{ccl}
H^{*,*}(\tilde D)&\stackrel{h_{*}}{\longrightarrow}& H^{*,*}(\tilde {D'})\\
\vcenter{%
\llap {${t_{*}}$}}{{\Big\uparrow}}\;\;\Big\downarrow\vcenter{%
\rlap {$\pi_{*}$}}&& \Big\downarrow\vcenter{%
\llap {$\pi'_{*}\;\;$}}\;\vcenter{%
\rlap {$\;\;{t'_{*}}$}}{{\Big\uparrow}}\\
H^{*,*}_{G}(D)&\stackrel{\overline{h}_{*}}
{\longrightarrow}&H_{G}^{*,*}(D')
\end{array}}$
\end{center}
where $t_*$ (respectively $t'_*$) stands for the transfer map
corresponding to the action of $G$ on $\tilde D$ (respectively
$\tilde {D'}$). Since we are working with coefficients in $\F_2$ and
the order of $G$ is odd, then both $\pi_*t_*$ and  $\pi'_*t'_*$ are
isomorphisms. In addition, the commutative diagram implies that
$h_*t_*=t'_*\overline h_*$. Using the fact that $h_*$ is an
isomorphism we should be able to conclude that $\overline h_*$ is
injective. A similar argument using the  fact that $\pi_*$ and
$\pi'_*$ are onto implies that $\overline h_*$ is surjective.
Finally, the equivariant homologies
of $D$ and $D'$ are isomorphic. \\
The invariance under the left twisted  first Reidemeister   move is
proved in a similar way.
\subsubsection*{5.2- Invariance under the second Reidemeister move} We will
switch to framed links for a while. Let $D$ and $D'$ be two link
diagrams related by a single second Reidemeister move and assume
that $D'$ is the one that has more crossings, see figure 5. Let
$\tilde D$ and $\tilde D'$ be the two covering diagrams. We shall
prove that the equivariant homologies are isomorphic. Here we
consider the framed version $H_{*,*}^G (D)$ and $H_{*,*}^G ( D')$.
\begin{center}
\includegraphics[width=4cm,height=1cm]{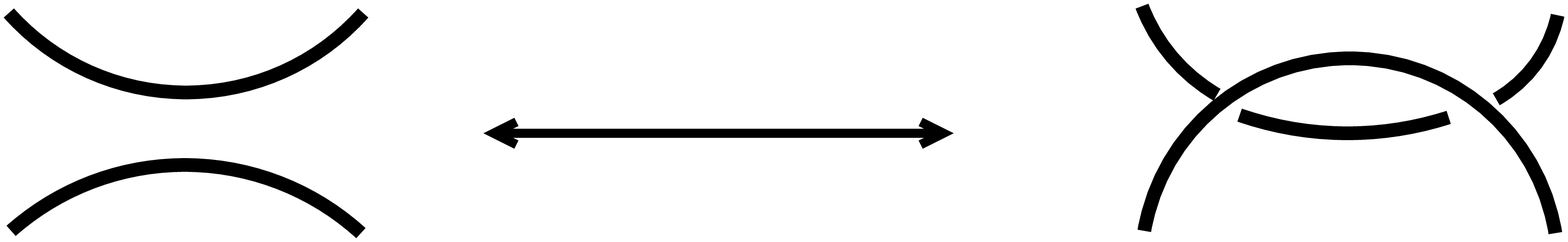}\\
Figure 5.
\end{center}
Following \cite{APS}, we define two maps $\overline \alpha$ and
$\overline \beta$  $$\func{\overline
\beta}{C_{p,q-2}(\includegraphics[width=0.5cm,height=0.5cm]{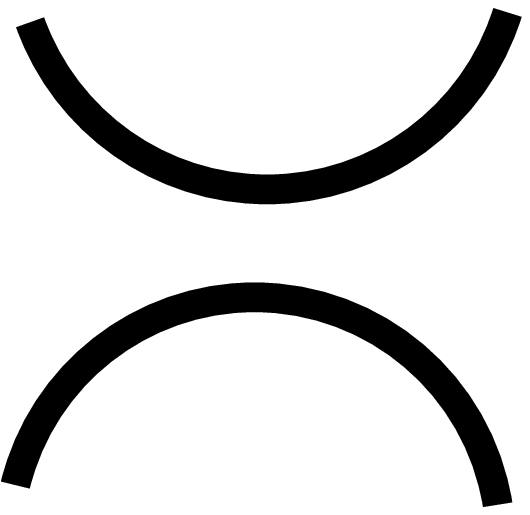})}{C_{p,q-1}(\includegraphics[width=0.5cm,height=0.5cm]{kplus})}{
\includegraphics[width=0.5cm,height=0.5cm]{kauffzero}}
{\includegraphics[width=0.5cm,height=0.5cm]{positivemarker}} $$ and

$$\doublefunc{\overline
\alpha}{C_{p,q-1}(\includegraphics[width=0.5cm,height=0.5cm]{kplus})}{C_{p,q}(\includegraphics[width=0.5cm,height=0.5cm]{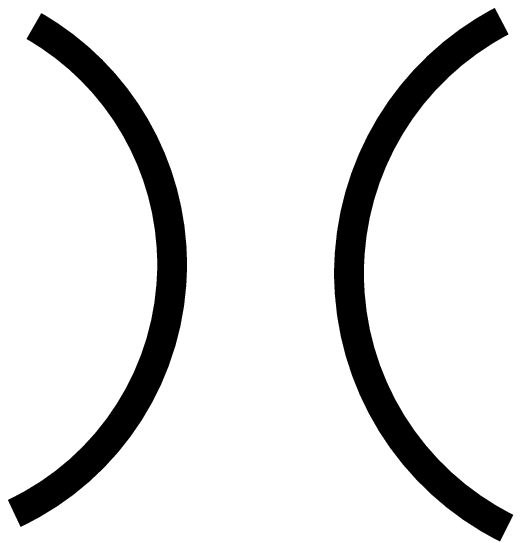})}{\includegraphics[width=0.5cm,height=0.5cm]{positivemarker}}{0}
{\includegraphics[width=0.5cm,height=0.5cm]{negativemarker}}{\includegraphics[width=0.5cm,height=0.5cm]{kauffinfini}}$$
Now, we set: $\gamma=\overline \alpha d_v \overline \beta$ which is
chain map from
$C_{p,q}(\includegraphics[width=0.5cm,height=0.5cm]{kauffzero})\longmapsto
C_{p-1,q+2}(\includegraphics[width=0.5cm,height=0.5cm]{kauffinfini})$.
We define two maps $f$ and $g$ as follows:
$$\func{f}{C_{p,q}(\includegraphics[width=0.5cm,height=0.5cm]{kauffzero})}{C_{p,q}(D')}{\includegraphics[width=1cm,height=1cm]{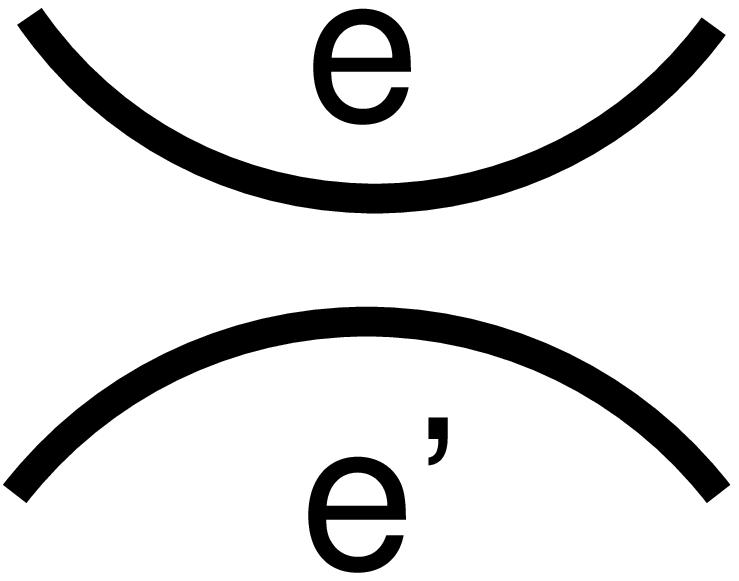}}
{\includegraphics[width=1cm,height=1cm]{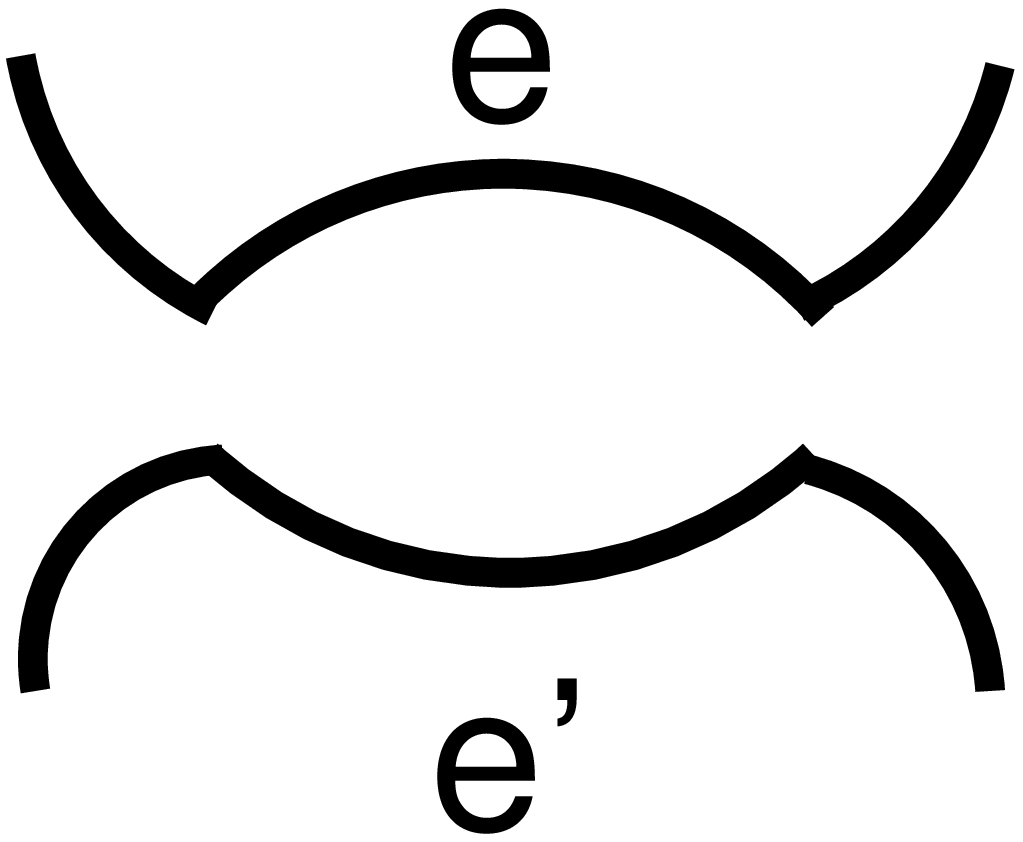}}$$ and
$$\func{g}{C_{p,q}(\includegraphics[width=0.5cm,height=0.5cm]{kauffinfini})}{C_{p+1,q-2}(D')}
{\includegraphics[width=1cm,height=1cm]{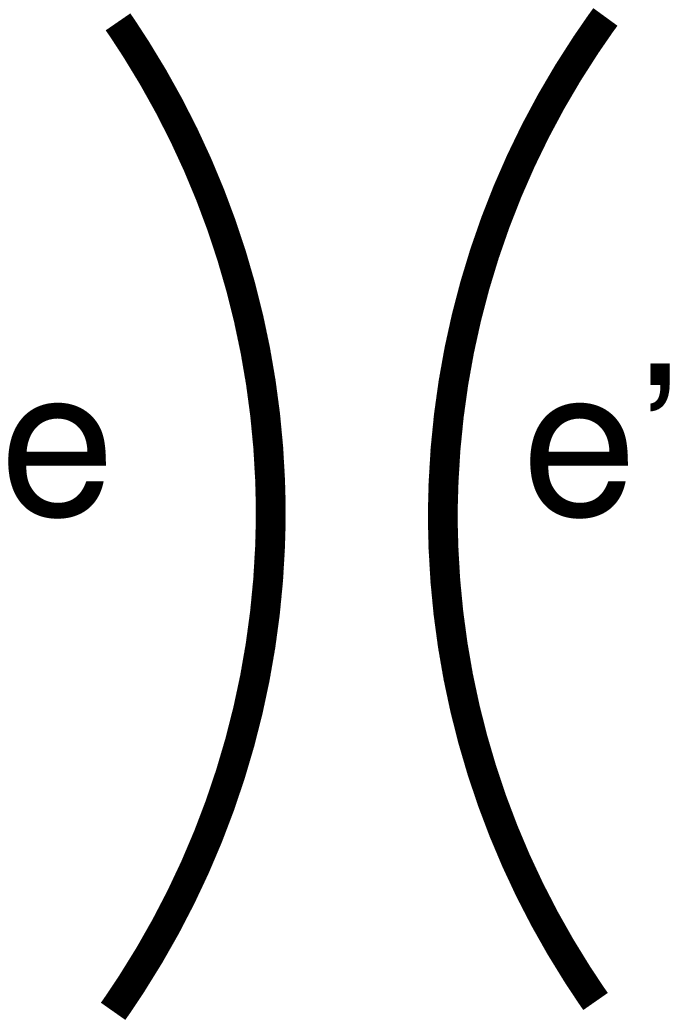}}{\includegraphics[width=1cm,height=1cm]{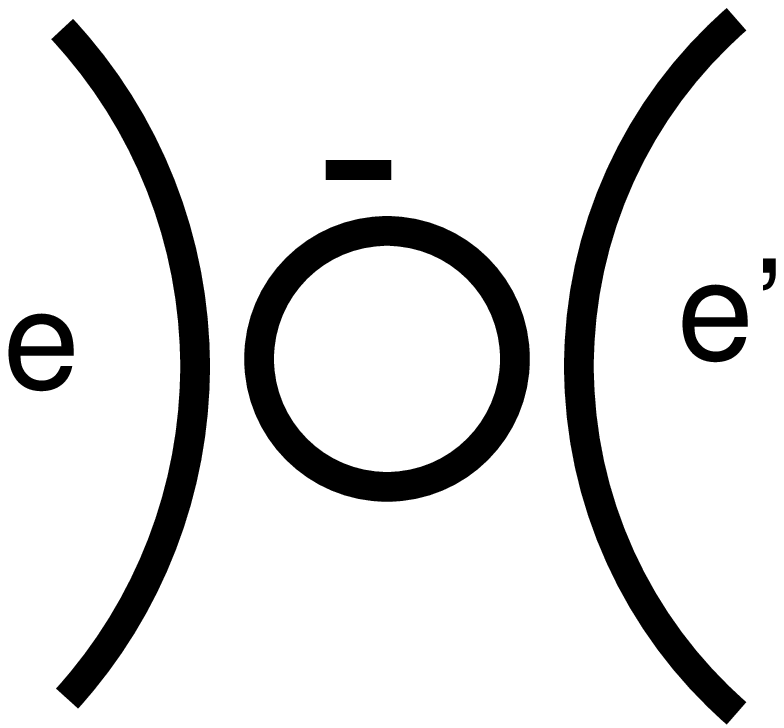}.}$$
Let $\rho$ be the chain  map $\rho=f+g\circ\gamma:
C_{p,q}(D)\longmapsto
C_{p,q}(D')$.\\
\textbf{Theorem 5.2. \cite{APS}} \emph{ The map $\rho$ induces an
isomorphism
in homology.}\\
The diagrams $\tilde D$ and  $\tilde D'$ differ  by $p$ Reidemeister
moves of type 2. To each move we associate a map $\rho$ as explained
earlier. Let us label these maps by
$\rho_v,\rho_{\varphi(v)},\dots,\rho_{\varphi^{p-1}(v)}$. By
composing these maps we define a map
$\Phi=\rho_v\circ\rho_{\varphi(v)}\circ\dots\circ\rho_{\varphi^{p-1}(v)}:
C_{p,q}(\tilde D)\longmapsto C_{p,q}(\tilde {D'})$. It is easy to
see that we have $\rho_w\circ \rho_{w'} =\rho_{w'}\circ \rho_w$ and
$\varphi \circ \rho_w=\rho_{\varphi(w)}\circ \varphi$. Consequently,
$\Phi$ is $G-$equivariant. Thus, it induces a map $\overline \Phi$
between the quotient chain complexes. Arguments similar to those
used in the case of the invariance under first Reidemeister move
should enable us to conclude that $\overline \Phi$ induces an
isomorphism between the equivariant Khovanov homologies of $D$ and
$D'$. \subsubsection*{5.3- Invariance under the third Reidemeister
move} In this paragraph, we shall prove the invariance of our
equivariant homology under the third Reidemeister move. Once again,
we are going to work with the framed version. Let $D$ and $D'$ be
two diagrams which differ by a single third Reidemeister move as in
the following picture
\begin{center}
\includegraphics[width=3cm,height=1cm]{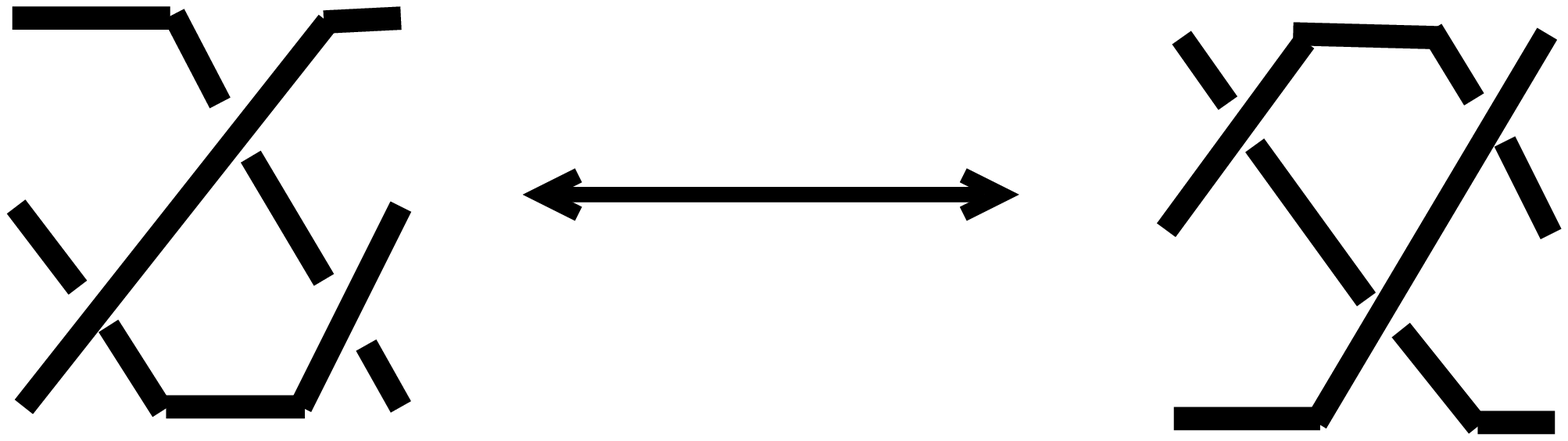}\\
Figure 6.
\end{center}
Our proof is based both on the construction in \cite{APS} and the
techniques we have developed in the previous paragraph. It could be
helpful if the reader has a copy of \cite{APS} with him. Let us
consider the following diagrams $D_{+}=$
\includegraphics[width=0.75cm,height=0.5cm]{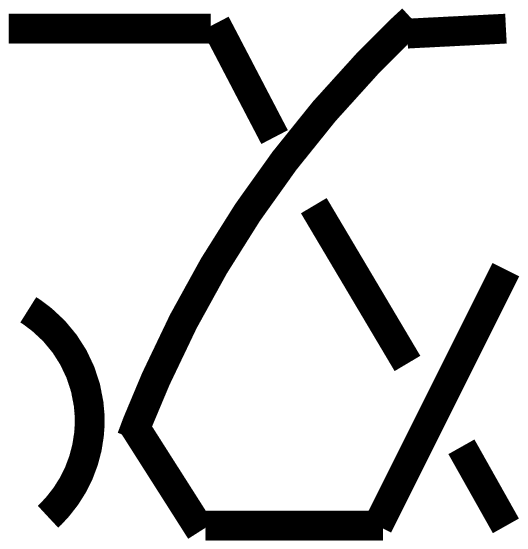}, $D_{-}=$ \includegraphics[width=0.75cm,height=0.5cm]{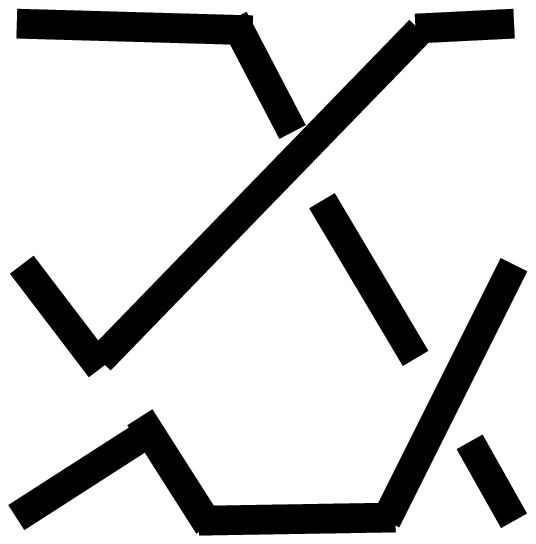} and $D_{++-}=$
\includegraphics[width=0.75cm,height=0.5cm]{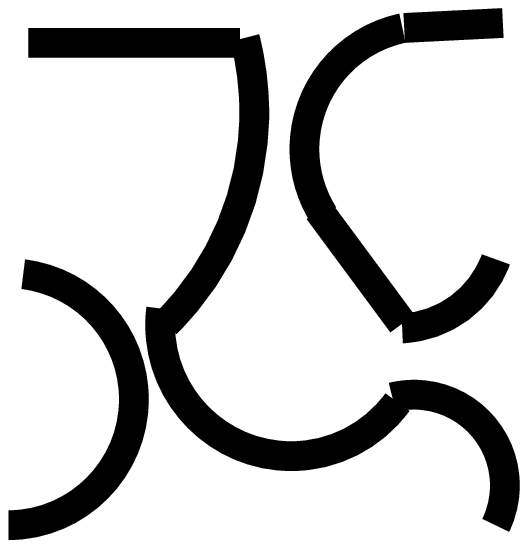}. We define
diagrams  $D'_{+}$, $D'_{-}$ and $D'_{++-}$ in the same way. Note
that the signs in the subscripts refer to the marker associated to
the
considered crossing.\\
The diagrams $D_{++-}$ and $D_{+}$ differ by a single Reidemeister
move of type 2. As we have explained in the previous paragraph there
exists a map $\Phi: C_{*,*}(\tilde D_{++-} )\longmapsto
C_{*,*}(\tilde D_{+} ) $ which is $G-$equivariant and induces an
isomorphism in homology. Now, let $C'_{*,*}(\tilde
{D_{+}})=\Phi(C_{*,*}(\tilde D_{++-}) )$ and consider the map
$\tilde i: C'_{*,*}(\tilde {D_{+}})\longmapsto C_{*,*}(\tilde
{D_{+}}) $.

We set $\tilde \beta: C_{*,*}(\tilde D) \longmapsto C_{*,*}(\tilde
D_+)$ to be the map define by:
\begin{center}
$\includegraphics[width=0.5cm,height=0.5cm]{positivemarker}
\dots\includegraphics[width=0.5cm,height=0.5cm]{positivemarker}
\longmapsto
\includegraphics[width=0.5cm,height=0.5cm]{kauffzero}\dots \includegraphics[width=0.5cm,height=0.5cm]{kauffzero}$
\end{center} and zero otherwise. Let $C_{*,*}'(\tilde D)=\tilde \beta^{-1}(C'_{*,*}(\tilde
D_+))$. We have the following \cite{APS}\\

\textbf{Lemma 5.3.} \emph{The maps $\tilde i: C'_{*,*}(\tilde D_+)
\longmapsto C_{*,*}(\tilde D_+)
 $ (resp. $\tilde {i'}: C'_{*,*}(\tilde D'_+) \longmapsto C_{*,*}(\tilde D'_+)
 $  ) and $\tilde j: C'_{*,*}(\tilde D) \longmapsto C_{*,*}(\tilde
 D)$ (resp. $\tilde {j'}: C'_{*,*}(\tilde {D'}) \longmapsto C_{*,*}(\tilde {D'})
 $)
  induce isomorphisms in homology. }\\
 \emph{Proof.} The induced map $\tilde{i}_*$ is an isomorphism in
 homology because it is a composition of isomorphisms, see \cite[Proposition 11.10]{APS}. Same argument applies for $\tilde{j}_*$. \fin\\

 Similarly to the case of the second Reidemeister move discussed earlier,
  by composition of the maps of type $\rho_{III}$ defined in
 \cite{APS} we should be able to construct a map $\tilde \Psi: C'_{*,*}(\tilde D) \longmapsto
 C'_{*,*}(\tilde D')$ which is $G-$equivariant and induces an isomorphism
 in homology. This map induces an isomorphism in homology $\overline
 \Psi$ between the homology of $\overline {C'_{*,*}(\tilde D)}$ and the
 homology of $\overline {C'_{*,*}(\tilde D')}$. Consequently, $\overline {\tilde {j'}_*}\circ\overline
 {\Psi_*}
 \circ \overline{\tilde {j}_*^{-1}}: H^G_{*,*}(D)\longmapsto H^G_{*,*}(D')$ is an
isomorphism. This completes the proof of the invariance under the
third  Reidemeister move.\\
Finally we use Theorem 4.2 to prove that  $H_G^{*,*}(L)$ is
isomorphic to $H^{*,*}(\tilde L)^G$. This completes the proof of
Theorem 1.
\subsection*{6- Equivariant Khovanov homology for framed links in the solid
torus} In this section, we show how our equivariant construction can
be described in the context of the categorification of the Kauffman
bracket skein module of the solid torus \cite{APS}. Everything here
is done similarly to what we have discussed in the previous
sections. Consequently, we
are going to  omit the details and describe things briefly. We first review the notion of skein modules.\\
Let $M$ be an oriented compact three-manifold. A framed link in $M$
is an embedding of a finite family of annuli into the interior of
$M$. Let $\cal L$ be the set of all isotopy classes of framed
 links in $M$ including the
empty link. Let  $\Z[A^{\pm}]$$ {\cal L}$ be the free module
generated by $\cal
 L$. The Kauffman bracket skein module of $M$, denoted here by $\cal K$$(M)$, is defined as the quotient of
$\Z[A^{\pm}]$${\cal L}$ by the smallest submodule generated by all elements of the following form\\
1) $L\cup \bigcirc + (A^2+A^{-2}) L $, where $L$ is any framed link
in $M$, and
$L \cup \bigcirc$ is the disjoint union of $L$ with a trivial component,\\
2) $L_+ - AL_0-A^{-1}L_{\infty}$,
 where $L_+$, $L_0$ and $L_{\infty}$ are three links which are
 identical except in a small ball where they look like in figure
 3.\\
 The existence and the uniqueness of the Jones polynomial is
 equivalent to the fact that ${\cal{K}}(S^3)$ is isomorphic to
 $\Z[A^{\pm}]$ with the empty link as a basis. If $F$ is an oriented
 surface, then the skein module of $F\times I$ admits an algebra structure \cite{Bu}. In particular,
 the   skein algebra of the solid torus
$S^1 \times I \times I$ is isomorphic to the polynomial algebra
$\Z[A^{\pm}][z]$, where $z$ is represented by a nontrivial curve in
the annulus as in the following picture \null
\begin{center}
\begin{picture}(0,0)
 \put(0,0){\circle{10}} \put(0,0){\circle{20}} \put(12,0){$z$}
 \put(0,0){\circle{40}}
\end{picture}
\end{center}
\vspace{0.2cm}

Let $L$ be a link in the solid torus. Let  $D$ be  diagram of $L$ in
the annulus and let $(C^{*,*,*}(D),d)$ be the chain complex of $D$
with coefficients in $\F_2$. The skein module of the solid torus has
a basis made up of links of the form $z^n$ ($n$ parallel copies of
$z$). Thus, we shall use $n$ as the third script instead of $z^n$ as
in  original definition \cite{APS}. The homology of
$(C^{*,*,*}(D),\F_2)$ defines an invariant of framed links in the
solid
torus.\\
Let $\tilde L$ be the covering link of $L$ in the $p-$fold cyclic
cover of the solid torus.  Let $\tilde D$ be a symmetric diagram of
$\tilde L$. Arguments similar to those used in section 3 show that
the action of the rotation  on the diagram $\tilde D$ extends to an
action of the finite cyclic group $G$ on the chain complex
$(C^{*,*,*}(\tilde D),d)$. The homology $H^{*,*,*}_G$ of the
quotient complex $(\overline {C^{*,*,*}(\tilde D)},\overline{d})$ is
called the $G-$equivariant Khovanov homology of $D$. The proofs in
the previous section extend straightforward to conclude that
$H^{*,*,*}_G$ is an invariant of framed links in the solid torus.
\subsection*{7- Equivariant  graph homology }
In this section, we explain how one may extend our link equivariant
homology to graphs. Let us first fix notations and review some
definitions. Throughout the rest of this paper, a  graph is a
1-dimensional finite CW-complex. Let $\mathcal G$ be a graph with
vertex set $V({\cal G})$ and edge set $E({\mathcal G})$. The
chromatic polynomial of $\cal G$ is a one variable  polynomial
$P({\mathcal G})\in \Z[\lambda]$ which when evaluated at an integer
$m$ gives the number of colorings of the vertices of $\mathcal G$ by
a palette of $m-$colors satisfying the property  that  vertices
which are connected by a edge  have different colors. Now, we shall
briefly review the definition of graph homology following \cite{HR}.
We consider homology with coefficients in $\F_2$. Take a set of
colors $\{1,x\}$ and define a product $\star$ as in $\Z[x]/x^2$. For
each $s\subseteq E({\mathcal G})$, we set $[G:s]$ to be the graph
whose vertex set is $V({\cal G})$ and whose edge set is $s$. An
\emph{enhanced state} of $\mathcal {G}$ is $S=(s,c)$ where
$s\subseteq E({\mathcal G})$ and $c$ is an assignment of $1$ or $x$
to each connected component of the spanning subgraph $[G:s]$. If $S$
is an enhanced state then we set $i(S)$ to be the number of edges in
$S$ and we set $j(S)$ to be the number of $x$'s in $c$. Now, we
define $C^{i,j}({\mathcal G})$ to be the vector space generated by
all enhanced states of $\mathcal G$ with $i(S)=i$ and $j(S)=j$. The
differential is defined by
$$\func{d} {C^{i,j}({\mathcal G})} {C^{i+1,j}({\mathcal G})} {S}
{\displaystyle\sum_{e \in E({\mathcal G}-s)}S_e,}$$
 where $S_e$ is any enhanced state obtained from $S$ by adding an
 edge not in $s$ and adjusting the sign according to the product
 $\star$,  see \cite[Page 1375]{HR} for more
 details. The homology of $C^{*,*}({\mathcal G},d)$ is an invariant
 of $\mathcal G$. The chromatic polynomial  is the Euler
 characteristic of $H^{*,*}({\mathcal G})$ evaluated at
 $q=\lambda-1$.\\
 Let $\tilde{\mathcal G}$ be a graph on which  the finite cyclic
 group
 $G$ acts. The action of $G$
 on $\tilde{\mathcal G}$ extends to an action on the set of enhanced states.
 Thus, the group $G$ acts on $C^{i,j}(\tilde {\mathcal G})$. The action of
 $G$ commutes with the differential for the same reason as in the
 proof of Lemma 3.1. We obtain a quotient complex, its homology is
 called the $G-$equivariant homology of $\tilde {\mathcal G}$.

\end{document}